# Third-order Smoothness Helps: Even Faster Stochastic Optimization Algorithms for Finding Local Minima

Yaodong Yu[*][‡] and Pan Xu[†][‡] and Quanquan Gu[§]

## Abstract

We propose stochastic optimization algorithms that can find local minima faster than existing algorithms for nonconvex optimization problems, by exploiting the third-order smoothness to escape non-degenerate saddle points more efficiently. More specifically, the proposed algorithm only needs $\widetilde{O}(\epsilon^{-10/3})$ stochastic gradient evaluations to converge to an approximate local minimum $\mathbf{x}$, which satisfies $\|\nabla f(\mathbf{x})\|_2 \leq \epsilon$ and $\lambda_{\min}(\nabla^2 f(\mathbf{x})) \geq -\sqrt{\epsilon}$ in the general stochastic optimization setting, where $\widetilde{O}(\cdot)$ hides logarithm polynomial terms and constants. This improves upon the $\widetilde{O}(\epsilon^{-7/2})$ gradient complexity achieved by the state-of-the-art stochastic local minima finding algorithms by a factor of $\widetilde{O}(\epsilon^{-1/6})$. For nonconvex finite-sum optimization, our algorithm also outperforms the best known algorithms in a certain regime.

## 1 Introduction

We study the following problem of unconstrained stochastic optimization

$$\min_{\mathbf{x}\in\mathbb{R}^d} f(\mathbf{x}) = \mathbb{E}_{\xi\sim\mathcal{D}}[F(\mathbf{x};\xi)], \tag{1.1}$$

where $F(\mathbf{x};\xi):\mathbb{R}^d \to \mathbb{R}$ is a stochastic function and $\xi$ is a random variable sampled from a fixed distribution $\mathcal{D}$. In particular, we are interested in nonconvex optimization where the expected function $f(\mathbf{x})$ is nonconvex. Note that (1.1) also covers the finite sum optimization problem, where $f(\mathbf{x}) = 1/n\sum_{i=1}^n f_i(\mathbf{x})$, which can be seen as a special case of stochastic optimization when the random variable $\xi$ follows a uniform distribution over $\{1,\ldots,n\}$. This kind of nonconvex optimization is ubiquitous in machine learning, especially deep learning (LeCun et al., 2015). Finding a global minimum of nonconvex problem (1.1) is generally NP hard (Hillar and Lim, 2013). Nevertheless, for many nonconvex optimization in machine learning, a local minimum is adequate and can be as good as a global minimum in terms of generalization performance, especially in many deep learning problems (Choromanska et al., 2015; Dauphin et al., 2014).

---

[*]Department of Computer Science, University of Virginia, Charlottesville, VA 22904, USA; e-mail:yy8ms@virginia.edu

[†]Department of Computer Science, University of Virginia, Charlottesville, VA 22904, USA; e-mail: px3ds@virginia.edu

[‡]Equal Contribution

[§]Department of Computer Science, University of Virginia, Charlottesville, VA 22904, USA; e-mail: qg5w@virginia.edu

In this paper, we aim to design efficient stochastic optimization algorithms that can find an approximate local minimum of (1.1), i.e., an $(\epsilon, \epsilon_H)$-second-order stationary point defined as follows

$$\|\nabla f(\mathbf{x})\|_2 \leq \epsilon, \text{ and } \lambda_{\min}(\nabla^2 f(\mathbf{x})) \geq -\epsilon_H, \tag{1.2}$$

where $\epsilon, \epsilon_H \in (0, 1)$. Notably, when $\epsilon_H = \sqrt{L_2 \epsilon}$ for Hessian Lipschitz $f(\cdot)$ with parameter $L_2$, (1.2) is equivalent to the definition of $\epsilon$-second-order stationary point (Nesterov and Polyak, 2006). Algorithms based on cubic regularized Newton's method (Nesterov and Polyak, 2006) and its variants (Agarwal et al., 2016; Carmon and Duchi, 2016; Curtis et al., 2017; Kohler and Lucchi, 2017; Xu et al., 2017; Tripuraneni et al., 2017) have been proposed to find the approximate local minimum. However, all of them need to solve the cubic problems exactly (Nesterov and Polyak, 2006) or approximately (Agarwal et al., 2016; Carmon and Duchi, 2016) in each iteration, which poses a rather heavy computational overhead. Another line of research employs the negative curvature direction to find the local minimum by combining accelerated gradient descent and negative curvature descent (Carmon et al., 2016; Allen-Zhu, 2017), which becomes impractical in large scale and high dimensional machine learning problems due to the frequent computation of negative curvature in each iteration.

To alleviate the computational burden of local minimum finding algorithms, there has emerged a fresh line of research (Xu and Yang, 2017; Allen-Zhu and Li, 2017b; Jin et al., 2017b; Yu et al., 2017) that tries to achieve the iteration complexity as the state-of-the-art second-order methods, while only utilizing first-order oracles. The key observation is that first-order methods with noise injection (Ge et al., 2015; Jin et al., 2017a) are essentially an equivalent way to extract the negative curvature direction at saddle points (Xu and Yang, 2017; Allen-Zhu and Li, 2017b). With first-order information and the state-of-the-art Stochastically Controlled Stochastic Gradient (SCSG) method (Lei et al., 2017), the aforementioned methods (Xu and Yang, 2017; Allen-Zhu and Li, 2017b) converge to an $(\epsilon, \sqrt{\epsilon})$-second-order stationary point (an approximate local minimum) within $\widetilde{O}(\epsilon^{-7/2})$ stochastic gradient evaluations, where $\widetilde{O}(\cdot)$ hides logarithm polynomial factors and constants. In this work, motivated by Carmon et al. (2017) which employed the third-order smoothness of $f(\cdot)$ in deterministic nonconvex optimization and proposed an algorithm that converges to a first-order stationary point with $\widetilde{O}(\epsilon^{-5/3})$ full gradient evaluations, we explore the benefits of third-order smoothness in finding an approximate local minimum in the stochastic nonconvex optimization setting. In particular, we propose a stochastic optimization algorithm that only utilizes first-order oracles and finds the $(\epsilon, \epsilon_H)$-second-order stationary point within $\widetilde{O}(\epsilon^{-10/3})$ stochastic gradient evaluations. Note that our gradient complexity matches that of the state-of-the-art stochastic optimization algorithm SCSG (Lei et al., 2017) for finding first-order stationary points. At the core of our algorithm is an exploitation of the third-order smoothness of the objective function $f(\cdot)$ which enables us to choose a larger step size in the negative curvature descent stage, and therefore leads to a faster convergence rate. The main contributions of our work are summarized as follows

- We show that third-order smoothness of nonconvex function can lead to a faster escape from saddle points in the stochastic optimization. We characterize, for the first time, the improvement brought by third-order smoothness in finding the negative curvature direction.

- We propose efficient stochastic algorithms for both finite-sum and general stochastic objective functions and prove faster convergence rates for finding local minima. More specifically, for finite-sum setting, the proposed algorithm converges to an $(\epsilon, \epsilon_H)$-second-order stationary point within $\widetilde{O}(n^{2/3}\epsilon^{-2} + n\epsilon_H^{-2} + n^{3/4}\epsilon_H^{-5/2})$ stochastic gradient evaluations, and for general stochastic optimization setting, our algorithm converges to an approximate local minimum with only $\widetilde{O}(\epsilon^{-10/3})$ stochastic gradient evaluations.

- In each outer iteration, our proposed algorithms only need to call either one negative curvature descent step, or an epoch of SCSG. Compared with existing algorithms which need to execute at least one of these two procedures or even both in each outer iteration, our algorithms save a lot of gradient and negative curvature computations.

The remainder of this paper is organized as follows: We review the related work in Section 2, and present some preliminary definitions in Section 3. In Section 4, we elaborate the motivation of using third-order smoothness and present faster negative curvature descent algorithms in finite-sum and stochastic settings. In Section 5, we present local minima finding algorithms and the corresponding theoretical analyses for both finite-sum and general stochastic nonconvex optimization problems. Finally, we conclude our paper in Section 6.

**Notation** We use $\mathbf{A} = [A_{ij}] \in \mathbb{R}^{d \times d}$ to denote a matrix and $\mathbf{x} = (x_1, ..., x_d)^\top \in \mathbb{R}^d$ to denote a vector. Let $\|\mathbf{x}\|_q = (\sum_{i=1}^d |x_i|^q)^{1/q}$ be $\ell_q$ vector norm for $0 < q < +\infty$. We use $\|\mathbf{A}\|_2$ and $\|\mathbf{A}\|_F$ to denote the spectral and Frobenius norm of $\mathbf{A}$. For a three-way tensor $\mathcal{T} \in \mathbb{R}^{d \times d \times d}$ and vector $\mathbf{x} \in \mathbb{R}^d$, we denote the inner product as $\langle \mathcal{T}, \mathbf{x}^{\otimes 3} \rangle$. For a symmetric matrix $\mathbf{A}$, let $\lambda_{\max}(\mathbf{A})$ and $\lambda_{\min}(\mathbf{A})$ be the maximum, minimum eigenvalues of matrix $\mathbf{A}$. We use $\mathbf{A} \succeq 0$ to denote $\mathbf{A}$ is positive semidefinite. For any two sequences $\{a_n\}$ and $\{b_n\}$, if $a_n \leq C b_n$ for some $0 < C < +\infty$ independent of $n$, we write $a_n = O(b_n)$. The notation $\widetilde{O}(\cdot)$ hides logarithmic factors. For sequences $f_n$, $g_n$, if $f_n$ is less than (larger than) $g_n$ up to a constant, then we write $f_n \lesssim g_n$ ($f_n \gtrsim g_n$).

## 2 Related Work

In this section, we discuss related work for finding approximate second-order stationary points in nonconvex optimization. Roughly speaking, existing literature can be divided into the following three categories.

**Hessian-based:** The pioneer work of Nesterov and Polyak (2006) proposed the cubic regularized Newton's method to find an $(\epsilon, \epsilon_H)$-second-order stationary point in $O(\max\{\epsilon^{-3/2}, \epsilon_H^{-3}\})$ iterations. Curtis et al. (2017) showed that the trust-region Newton method can achieve the same iteration complexity as the cubic regularization method. Recently, Kohler and Lucchi (2017); Xu et al. (2017) showed that by using subsampled Hessian matrix instead of the entire Hessian matrix in cubic regularization method and trust-region method, the iteration complexity can still match the original exact methods under certain conditions. However, these methods need to compute the Hessian matrix and solve a very expensive subproblem either exactly or approximately in each iteration, which can be computationally intractable for high-dimensional problems.

**Hessian-vector product-based:** Through different approaches, Carmon et al. (2016) and Agarwal et al. (2016) independently proposed algorithms that are able to find $(\epsilon, \sqrt{\epsilon})$-second-order stationary points with $\widetilde{O}(\epsilon^{-7/4})$ full gradient and Hessian-vector product evaluations. By making an additional assumption of the third-order smoothness on the objective function and combining the negative curvature descent with the "convex until proven guilty" algorithm, Carmon et al. (2017) proposed an algorithm that is able to find an $(\epsilon, \sqrt{\epsilon})$-second-order stationary point with $\widetilde{O}(\epsilon^{-5/3})$ full gradient and Hessian-vector product evaluations.[1] Agarwal et al. (2016) also considered the finite-sum nonconvex optimization and proposed algorithm which is able to find approximate local minima within $\widetilde{O}(n\epsilon^{-3/2} + n^{3/4}\epsilon^{-7/2})$ stochastic gradient and stochastic Hessian-vector product evaluations. For nonconvex finite-sum problems, Reddi et al. (2017) proposed an algorithm, which is a combination of first-order and second-order methods to find approximate $(\epsilon, \epsilon_H)$-second-order stationary points, and requires $\widetilde{O}(n^{2/3}\epsilon^{-2} + n\epsilon_H^{-3} + n^{3/4}\epsilon_H^{-7/2})$ stochastic gradient and stochastic Hessian-vector product evaluations. In the general stochastic optimization setting, Allen-Zhu (2017) proposed an algorithm named Natasha2, which is based on variance reduction and negative curvature descent, and is able to find $(\epsilon, \sqrt{\epsilon})$-second-order stationary points with at most $\widetilde{O}(\epsilon^{-7/2})$ stochastic gradient and stochastic Hessian-vector product evaluations. Tripuraneni et al. (2017) proposed a stochastic cubic regularization algorithm to find $(\epsilon, \sqrt{\epsilon})$-second-order stationary points and achieved the same runtime complexity as Allen-Zhu (2017).

**Gradient-based:** For general nonconvex problems, Levy (2016); Jin et al. (2017a,b) showed that it is possible to escape from saddle points and find local minima by using gradient evaluations plus random perturbation. The best-known runtime complexity of these methods is $\widetilde{O}(\epsilon^{-7/4})$ when setting $\epsilon_H = \sqrt{\epsilon}$, which is proposed in Jin et al. (2017b). In the nonconvex finite-sum setting, Allen-Zhu and Li (2017b) proposed a first-order negative curvature finding method called Neon2 and applied it to the nonconvex stochastic variance reduced gradient (SVRG) method (Reddi et al., 2016; Allen-Zhu and Hazan, 2016; Lei et al., 2017), which gives rise to an algorithm to find $(\epsilon, \epsilon_H)$-second-order stationary points with $\widetilde{O}(n^{2/3}\epsilon^{-2} + n\epsilon_H^{-3} + n^{3/4}\epsilon_H^{-7/2} + n^{5/12}\epsilon^{-2}\epsilon_H^{-1/2})$ stochastic gradient evaluations. For nonconvex stochastic optimization problems, a variant of stochastic gradient descent (SGD) in Ge et al. (2015) is proved to find the $(\epsilon, \sqrt{\epsilon})$-second-order stationary point with $O(\epsilon^{-4}\text{poly}(d))$ stochastic gradient evaluations. Recently, Xu and Yang (2017); Allen-Zhu and Li (2017b) turned the stochastically controlled stochastic gradient (SCSG) method (Lei et al., 2017) into approximate local minima finding algorithms, which involves stochastic gradient computation. The runtime complexity of these algorithms is $\widetilde{O}(\epsilon^{-10/3} + \epsilon^{-2}\epsilon_H^{-3})$. Very recently, in order to further save gradient and negative curvature computations, Yu et al. (2017) proposed a family of new algorithms called GOSE, which can find $(\epsilon, \epsilon_H)$-second-order stationary points with $\widetilde{O}(\epsilon^{-7/4} + (\epsilon^{-1/4} + \epsilon_H^{-1/2}) \min\{\epsilon_H^{-3}, N_\epsilon\})$ gradient evaluations in the deterministic setting, $\widetilde{O}(n^{2/3}\epsilon^{-2} + (n + n^{3/4}\epsilon_H^{-1/2}) \min\{\epsilon_H^{-3}, N_\epsilon\})$ stochastic gradient evaluations in the finite-sum setting, and $\widetilde{O}(\epsilon^{-10/3} + \epsilon^{-2}\epsilon_H^{-3} + (\epsilon^{-2} + \epsilon_H^{-2}) \min\{\epsilon_H^{-3}, N_\epsilon\})$ stochastic gradient evaluations in the general stochastic setting, where $N_\epsilon$ is the number of saddle points encountered by the algorithms until they find an approximate local minimum. Provided that $N_\epsilon$ is smaller than $\epsilon_H^{-3}$, the runtime complexity is better than the state-of-the-art (Allen-Zhu and Li, 2017b; Jin et al., 2017b). However, none of the above gradient-based algorithms explore the third-order smoothness of the nonconvex objective function.

---

[1] As shown in Carmon et al. (2017), the second-order accuracy parameter $\epsilon_H$ can be set as $\epsilon^{2/3}$ and the total runtime complexity remains the same, i.e., $\widetilde{O}(\epsilon^{-5/3})$.

# 3 Preliminaries

In this section, we present some definitions which will be used in our algorithm design and later theoretical analysis.

**Definition 3.1** (Smoothness). A differentiable function $f(\cdot)$ is $L_1$-smooth, if for any $\mathbf{x}, \mathbf{y} \in \mathbb{R}^d$:
$$\|\nabla f(\mathbf{x}) - \nabla f(\mathbf{y})\|_2 \leq L_1 \|\mathbf{x} - \mathbf{y}\|_2.$$

**Definition 3.2** (Hessian Lipschitz). A twice-differentiable function $f(\cdot)$ is $L_2$-Hessian Lipschitz, if for any $\mathbf{x}, \mathbf{y} \in \mathbb{R}^d$:
$$\|\nabla^2 f(\mathbf{x}) - \nabla^2 f(\mathbf{y})\|_2 \leq L_2 \|\mathbf{x} - \mathbf{y}\|_2.$$

Note that the Hessian-Lipschitz is also referred to as second-order smoothness. The above two smoothness conditions are widely used in nonconvex optimization problems (Nesterov and Polyak, 2006). In this paper, we will explore the third-order derivative Lipschitz condition. Here we use a three-way tensor $\nabla^3 f(\mathbf{x}) \in \mathbb{R}^{d \times d \times d}$ to denote the third-order derivative of a function, which is formally defined below.

**Definition 3.3** (Third-order Derivative). The third-order derivative of function $f(\cdot)$: $\mathbb{R}^d \to \mathbb{R}$ is a three-way tensor $\nabla^3 f(\mathbf{x}) \in \mathbb{R}^{d \times d \times d}$ which is defined as
$$[\nabla^3 f(\mathbf{x})]_{ijk} = \frac{\partial}{\partial x_i \partial x_j \partial x_k} f(\mathbf{x}) \quad \text{for } i,j,k = 1, \ldots, d \text{ and } \mathbf{x} \in \mathbb{R}^d.$$

Next we introduce the third-order smooth condition of function $f(\cdot)$, which implies that the third-order derivative will not change rapidly.

**Definition 3.4** (Third-order derivative Lipschitz). A thrice-differentiable function $f(\cdot)$ has $L_3$-Lipschitz third-order derivative, if for any $\mathbf{x}, \mathbf{y} \in \mathbb{R}^d$:
$$\|\nabla^3 f(\mathbf{x}) - \nabla^3 f(\mathbf{y})\|_F \leq L_3 \|\mathbf{x} - \mathbf{y}\|_2.$$

The above definition has been introduced in Anandkumar and Ge (2016), and the third-order derivative Lipschitz is also referred to as third-order smoothness in Carmon et al. (2017). One can also use another equivalent notion of third-order derivative Lipschitz condition used in Carmon et al. (2017). Note that the third-order Lipschitz condition is critical in our algorithms and theoretical analysis in later sections. In the sequel, we will use third-order derivative Lipschitz and third-order smoothness interchangeably.

**Definition 3.5** (Optimal Gap). For a function $f(\cdot)$, we define $\Delta_f$ as
$$f(\mathbf{x}_0) - \inf_{\mathbf{x} \in \mathbb{R}^d} f(\mathbf{x}) \leq \Delta_f,$$
and without loss of generality, we throughout assume $\Delta_f < +\infty$.

**Definition 3.6** (Geometric Distribution). For a random integer $X$, define $X$ has a geometric distribution with parameter $p$, denoted as $\text{Geom}(p)$, if it satisfies that
$$\mathbb{P}(X = k) = p^k(1-p), \quad \forall k = 0, 1, \ldots.$$

**Definition 3.7** (Sub-Gaussian Stochastic Gradient)**.** For any $\mathbf{x} \in \mathbb{R}^d$ and random variable $\xi \in \mathcal{D}$, the stochastic gradient $\nabla F(\mathbf{x}; \xi)$ is sub-Gaussian with parameter $\sigma$ if it satisfies

$$\mathbb{E}\bigg[\exp\bigg(\frac{\|\nabla F(\mathbf{x}; \xi) - \nabla f(\mathbf{x})\|_2^2}{\sigma^2}\bigg)\bigg] \leq \exp(1).$$

In addition, we introduce $\mathbb{T}_g$ to denote the time complexity of stochastic gradient evaluation, and $\mathbb{T}_h$ to denote the time complexity of stochastic Hessian-vector product evaluation.

## 4 Exploiting Third-order Smoothness

In this section we will show how to exploit the third-order smoothness of the objective function to make better use of the negative curvature direction for escaping saddle points. To begin with, we will first explain why third-order smoothness helps in general nonconvex optimization problems. Then we present two algorithms which are able to utilize the third-order smoothness to take a larger step size in the finite-sum setting and the stochastic setting.

In order to find local minima in nonconvex problems, different kinds of approaches have been explored to escape from saddle points. One of these approaches is to use negative curvature direction (Moré and Sorensen, 1979) to escape from saddle points, which has been explored in many previous studies (Carmon et al., 2016; Allen-Zhu, 2017). According to some recent work (Xu and Yang, 2017; Allen-Zhu and Li, 2017b), it is possible to extract the negative curvature direction by only using (stochastic) gradient evaluations, which makes the negative curvature descent approach more appealing.

We first consider a simple case to illustrate how to utilize the third-order smoothness when taking a negative curvature descent step. For nonconvex optimization problems, an $\epsilon$-first-order stationary point $\widehat{\mathbf{x}}$ can be found by using first-order methods, such as gradient descent. Suppose an $\epsilon$-first-order stationary point $\widehat{\mathbf{x}}$ is not an $(\epsilon, \epsilon_H)$-first-order stationary point, which means that there exists a unit vector $\widehat{\mathbf{v}}$ such that

$$\widehat{\mathbf{v}}^\top \nabla^2 f(\widehat{\mathbf{x}})\, \widehat{\mathbf{v}} \leq -\frac{\epsilon_H}{2}.$$

Following the previous work (Carmon et al., 2016; Xu and Yang, 2017; Allen-Zhu and Li, 2017b; Yu et al., 2017), one can take a negative curvature descent step along $\widehat{\mathbf{v}}$ to escape from the saddle point $\widehat{\mathbf{x}}$, i.e.,

$$\widetilde{\mathbf{y}} = \underset{\mathbf{y} \in \{\mathbf{u}, \mathbf{w}\}}{\operatorname{argmin}} f(\mathbf{y}), \quad \text{where } \mathbf{u} = \widehat{\mathbf{x}} - \alpha\, \widehat{\mathbf{v}},\ \mathbf{w} = \widehat{\mathbf{x}} + \alpha\, \widehat{\mathbf{v}} \text{ and } \alpha = O\bigg(\frac{\epsilon_H}{L_2}\bigg). \tag{4.1}$$

Suppose the function $f(\cdot)$ is $L_1$-smooth and $L_2$-Hessian Lipschitz, then the negative curvature descent step (4.1) is guaranteed to make the function value decrease,

$$f(\widetilde{\mathbf{y}}) - f(\widehat{\mathbf{x}}) = -O\bigg(\frac{\epsilon_H^3}{L_2^2}\bigg). \tag{4.2}$$

Inspired by the previous work (Carmon et al., 2017), incorporating with an additional assumption that the objective function has $L_3$-Lipschitz third-order derivatives (third-order smoothness) can help achieve better convergence rate for nonconvex problems. More specifically, we can adjust the update as follows,

$$\widehat{\mathbf{y}} = \operatorname*{argmin}_{\mathbf{y} \in \{\mathbf{u}, \mathbf{w}\}} f(\mathbf{y}), \quad \text{where } \mathbf{u} = \widehat{\mathbf{x}} - \eta \widehat{\mathbf{v}}, \mathbf{w} = \widehat{\mathbf{x}} + \eta \widehat{\mathbf{v}} \text{ and } \eta = O\left(\sqrt{\frac{\epsilon_H}{L_3}}\right). \tag{4.3}$$

Note that the adjusted step size $\eta$ can be much larger than the step size $\alpha$ in (4.1). Moreover, the adjusted negative curvature descent step (4.3) is guaranteed to decrease the function value by a larger decrement, i.e.,

$$f(\widehat{\mathbf{y}}) - f(\widehat{\mathbf{x}}) = -O\left(\frac{\epsilon_H^2}{L_3}\right). \tag{4.4}$$

Compared with (4.2), the decrement in (4.4) can be substantially larger. In other words, if we make the additional assumption of the third-order smoothness, the negative curvature descent with larger step size will make more progress toward decreasing the function value.

In the following, we will exploit the benefits of third-order smoothness in escaping from saddle points for stochastic nonconvex optimization.

## 4.1 Special Stochastic Optimization with Finite-Sum Structures

We first consider the following finite-sum problem, which is a special case of (1.1).

$$\min_{\mathbf{x} \in \mathbb{R}^d} f(\mathbf{x}) = \frac{1}{n} \sum_{i=1}^{n} f_i(\mathbf{x}), \tag{4.5}$$

where $f_i(\cdot)$ is possibly nonconvex. Algorithms in this setting have access to the information of each individual function $f_i(\cdot)$ and the entire function $f(\cdot)$. For the finite-sum structure, variance reduction-based methods (Reddi et al., 2016; Allen-Zhu and Hazan, 2016; Garber et al., 2016) can be applied to improve the gradient complexity of different nonconvex optimization algorithms.

We first present two different approaches to find the negative curvature direction in the finite-sum setting. The approximate leading eigenvector computation methods (Kuczyński and Woźniakowski, 1992; Garber et al., 2016) can be used to estimate the negative curvature direction, which are denoted by FastPCA methods in the sequel. To find the negative curvature direction of function $f(\cdot)$ at point $\mathbf{x}$, we can estimate the eigenvector with the smallest eigenvalue of a Hessian matrix $\mathbf{H} = \nabla^2 f(\mathbf{x})$. Since the function $f(\cdot)$ is $L_1$-smooth, one can instead estimate the eigenvector with the largest eigenvalue of the shifted matrix $\mathbf{M} = L_1 \cdot \mathbf{I} - \mathbf{H}$. Therefore, we can apply FastPCA for finite-sum problems to find the such eigenvector of the matrix $\mathbf{M}$, and we describe its result in the following lemma.

**Lemma 4.1.** (Garber et al., 2016) Let $f(\mathbf{x}) = 1/n \sum_{i=1}^{n} f_i(\mathbf{x})$ where each component function $f_i(\mathbf{x})$ is twice-differentiable and $L_1$-smooth. For any given point $\mathbf{x} \in \mathbb{R}^d$, if $\lambda_{\min}(\nabla^2 f(\mathbf{x})) \leq -\epsilon_H$, then with probability at least $1 - \delta$, Shifted-and-Inverted power method returns a unit vector $\widehat{\mathbf{v}}$ satisfying

$$\widehat{\mathbf{v}}^\top \nabla^2 f(\mathbf{x}) \widehat{\mathbf{v}} < -\frac{\epsilon_H}{2},$$

with $O\big((n + n^{3/4}\sqrt{L_1/\epsilon_H}) \log^3(d) \log(1/\delta)\big)$ stochastic Hessian-vector product evaluations.

Different from FastPCA-type methods, there is another type of approaches to compute the negative curvature direction without using the Hessian-vector product evaluation. In detail, Xu and Yang (2017); Allen-Zhu and Li (2017b) proposed algorithms that are able to extract the negative curvature direction based on random perturbation and (stochastic) gradient evaluation, and we denote these methods by Neon2 in the sequel. Here we adopt the Neon2$^{\text{svrg}}$ (Allen-Zhu and Li, 2017b) for the finite-sum setting and present its result in the following lemma.

**Lemma 4.2.** (Allen-Zhu and Li, 2017b) Let $f(\mathbf{x}) = 1/n \sum_{i=1}^n f_i(\mathbf{x})$ and each component function $f_i(\mathbf{x})$ is $L_1$-smooth and $L_2$-Hessian Lipschitz continuous. For any given point $\mathbf{x} \in \mathbb{R}^d$, with probability at least $1 - \delta$, Neon2$^{\text{svrg}}$ returns $\widehat{\mathbf{v}}$ satisfying one of the following conditions,

- $\widehat{\mathbf{v}} = \bot$, then $\lambda_{\min}(\nabla^2 f(\mathbf{x})) \geq -\epsilon_H$.
- $\widehat{\mathbf{v}} \neq \bot$, then $\widehat{\mathbf{v}}^\top \nabla^2 f(\mathbf{x}) \widehat{\mathbf{v}} \leq -\epsilon_H/2$ with $\|\mathbf{v}\|_2 = 1$.

The total number of stochastic gradient evaluations is $O\big((n + n^{3/4}\sqrt{L_1/\epsilon_H})\log^2(d/\delta)\big)$.

According to the above two lemmas, for finite-sum nonconvex optimization problems, there exists an algorithm, denoted by ApproxNC-FiniteSum($\cdot$), which uses stochastic gradient or stochastic Hessian-vector product to find the negative curvature direction. In detail, ApproxNC-FiniteSum($\cdot$) will return a unit vector $\widehat{\mathbf{v}}$ satisfying $\widehat{\mathbf{v}}^\top \nabla^2 f(\mathbf{x}) \widehat{\mathbf{v}} \leq -\epsilon_H/2$ if $\lambda_{\min}(\nabla^2 f(\mathbf{x})) < \epsilon_H$, otherwise the algorithm will return $\widehat{\mathbf{v}} = \bot$.

Based on the negative curvature finding algorithms, next we present the negative curvature descent algorithm in Algorithm 1, which is able to make larger decrease by taking advantage of the third-order smoothness.

---

**Algorithm 1** **N**egative **C**urvature **D**escent with **3**-order Smoothness - NCD3-FiniteSum ($f(\cdot)$, $\mathbf{x}$, $L_1$, $L_2$, $L_3$, $\delta$, $\epsilon_H$)

1: Set $\eta = \sqrt{3\epsilon_H/L_3}$
2: $\widehat{\mathbf{v}} \leftarrow$ ApproxNC-FiniteSum($f(\cdot), \mathbf{x}, L_1, L_2, \delta, \epsilon_H$)
3: **if** $\widehat{\mathbf{v}} \neq \bot$
4:     generate a Rademacher random variable $\zeta$
5:     $\widehat{\mathbf{y}} \leftarrow \mathbf{x} + \zeta \eta \widehat{\mathbf{v}}$
6:     **return** $\widehat{\mathbf{y}}$
7: **else**
8:     **return** $\bot$

---

As we can see from Algorithm 1, the step size is proportional to the square root of the accuracy parameter, i.e., $O(\sqrt{\epsilon_H})$. In sharp contrast, in the setting without third-order Lipschitz continuous assumption, the step size of the negative curvature descent step is much smaller, i.e., $O(\epsilon_H)$. As a consequence, Algorithm 1 will make greater progress if $\lambda_{\min}(\nabla^2 f(\mathbf{x})) < -\epsilon_H$, and we describe the theoretical result in the following lemma.

**Lemma 4.3.** Let $f(\mathbf{x}) = 1/n \sum_{i=1}^n f_i(\mathbf{x})$ and each component function $f_i(\cdot)$ is $L_1$-smooth, $L_2$-Hessian Lipschitz continuous, and the third derivative of $f(\cdot)$ is $L_3$-Lipschitz, suppose $\epsilon_H \in (0, 1)$ and set step size as $\eta = \sqrt{3\epsilon_H/L_3}$. If the input $\mathbf{x}$ of Algorithm 1 satisfies $\lambda_{\min}(\nabla^2 f(\mathbf{x})) < -\epsilon_H$, then with probability $1 - \delta$, Algorithm 1 will return $\widehat{\mathbf{y}}$ such that

$$\mathbb{E}_\zeta[f(\mathbf{x}) - f(\widehat{\mathbf{y}})] \geq \frac{3\epsilon_H^2}{8L_3},$$

where $\mathbb{E}_\zeta$ denotes the expectation over the Rademacher random variable $\zeta$. Furthermore, if we choose $\delta \leq \epsilon_H/(3\epsilon_H + 8L_2)$, it holds that

$$\mathbb{E}[f(\widehat{\mathbf{y}}) - f(\mathbf{x})] \leq -\frac{\epsilon_H^2}{8L_3},$$

where $\mathbb{E}$ is over all randomness of the algorithm, and the total runtime is $\widetilde{O}\big((n + n^{3/4}\sqrt{L_1/\epsilon_H}\,)\mathbb{T}_h\big)$ if ApproxNC-FiniteSum adopts Shifted-and-Inverted power method, and $\widetilde{O}\big((n + n^{3/4}\sqrt{L_1/\epsilon_H}\,)\mathbb{T}_g\big)$ if ApproxNC-FiniteSum adopts Neon2$^{\text{svrg}}$.

## 4.2 General Stochastic Optimization

Now we consider the general stochastic optimization problem in (1.1). In this setting, one cannot have access to the full gradient and Hessian information. Instead, only stochastic gradient and stochastic Hessian-vector product evaluations are accessible. As a result, we have to employ stochastic optimization methods to calculate the negative curvature direction. Similar to the finite-sum setting, there exist two different types of methods to calculate the negative curvature direction, and we first describe the result of online PCA method, i.e., Oja's algorithm (Allen-Zhu and Li, 2017a), in the following lemma.

**Lemma 4.4.** (Allen-Zhu and Li, 2017a) Let $f(\mathbf{x}) = \mathbb{E}_{\xi \sim \mathcal{D}}[F(\mathbf{x}; \xi)]$, where each component function $F(\mathbf{x}; \xi)$ is twice-differentiable and $L_1$-smooth. For any given point $\mathbf{x} \in \mathbb{R}^d$, if $\lambda_{\min}(\nabla^2 f(\mathbf{x})) \leq -\epsilon_H$, then with probability at least $1 - \delta$, Oja's algorithm returns a unit vector $\widehat{\mathbf{v}}$ satisfying

$$\widehat{\mathbf{v}}^\top \nabla^2 f(\mathbf{x})\,\widehat{\mathbf{v}} < -\frac{\epsilon_H}{2},$$

with $O\big((L_1^2/\epsilon_H^2)\log^2(d/\delta)\log(1/\delta)\big)$ stochastic Hessian-vector product evaluations.

The Oja's algorithm used in the above lemma can be seen as a stochastic variant of the FastPCA method introduced in previous subsection. Next we present another method called Neon2$^{\text{online}}$ (Allen-Zhu and Li, 2017b) to calculate the negative curvature direction in the stochastic setting.

**Lemma 4.5.** (Allen-Zhu and Li, 2017b) Let $f(\mathbf{x}) = \mathbb{E}_{\xi \sim \mathcal{D}}[F(\mathbf{x}; \xi)]$ where each component function $F(\mathbf{x}; \xi)$ is $L_1$-smooth and $L_2$-Hessian Lipschitz continuous. For any given point $\mathbf{x} \in \mathbb{R}^d$, with probability at least $1 - \delta$, Neon2$^{\text{online}}$ returns $\widehat{\mathbf{v}}$ satisfying one of the following conditions,

- $\widehat{\mathbf{v}} = \bot$, then $\lambda_{\min}(\nabla^2 f(\mathbf{x})) \geq -\epsilon_H$.
- $\widehat{\mathbf{v}} \neq \bot$, then $\widehat{\mathbf{v}}^\top \nabla^2 f(\mathbf{x})\,\widehat{\mathbf{v}} \leq -\epsilon_H/2$ with $\|\mathbf{v}\|_2 = 1$.

The total number of stochastic gradient evaluations is $O\big((L_1^2/\epsilon_H^2)\log^2(d/\delta)\big)$.

Based on Lemmas 4.4 and 4.5, there exists an algorithm, denoted by ApproxNC-Stochastic($\cdot$), which uses stochastic gradient evaluations or stochastic Hessian-vector product evaluations to find the negative curvature direction in the stochastic setting. Specifically, ApproxNC-Stochastic($\cdot$) returns a unit vector $\widehat{\mathbf{v}}$ that satisfies $\widehat{\mathbf{v}}^\top \nabla^2 f(\mathbf{x}) \widehat{\mathbf{v}} \leq -\epsilon_H/2$ provided $\lambda_{\min}(\nabla^2 f(\mathbf{x})) < \epsilon_H$, otherwise it will return $\widehat{\mathbf{v}} = \bot$. Next we present our algorithm in Algorithm 2 built on the negative curvature finding algorithms in the stochastic setting.

---

**Algorithm 2** NCD3-Stochastic ($f(\cdot)$, $\mathbf{x}$, $L_1$, $L_2$, $L_3$, $\delta$, $\epsilon_H$)

1: Set $\eta = \sqrt{3\epsilon_H/L_3}$
2: $\widehat{\mathbf{v}} \leftarrow$ ApproxNC-Stochastic($f(\cdot), \mathbf{x}, L_1, L_2, \delta, \epsilon_H$)
3: **if** $\widehat{\mathbf{v}} \neq \bot$
4:     generate a Rademacher random variable $\zeta$
5:     $\widehat{\mathbf{y}} \leftarrow \mathbf{x} + \zeta \eta \widehat{\mathbf{v}}$
6:     **return** $\widehat{\mathbf{y}}$
7: **else**
8:     **return** $\bot$

---

Algorithm 2 is analogous to Algorithm 1 in the finite-sum setting, with the same step size for the negative curvature descent step, i.e., $\eta = O(\sqrt{\epsilon_H})$. As a result, Algorithm 2 can also make greater progress when $\lambda_{\min}(\nabla^2 f(\mathbf{x})) < -\epsilon_H$, and we summarize its theoretical guarantee in the following lemma.

**Lemma 4.6.** Let $f(\mathbf{x}) = \mathbb{E}_{\xi \sim \mathcal{D}}[F(\mathbf{x};\xi)]$ and each stochastic function $F(\mathbf{x};\xi)$ is $L_1$-smooth, $L_2$-Hessian Lipschitz continuous, and the third derivative of $f(\mathbf{x})$ is $L_3$-Lipschitz. Suppose $\epsilon_H \in (0,1)$ and set step size as $\eta = \sqrt{3\epsilon_H/L_3}$. If the input $\mathbf{x}$ of Algorithm 2 satisfies $\lambda_{\min}(\nabla^2 f(\mathbf{x})) < -\epsilon_H$, then with probability $1 - \delta$, Algorithm 2 will return $\widehat{\mathbf{y}}$ such that

$$\mathbb{E}_\zeta[f(\mathbf{x}) - f(\widehat{\mathbf{y}})] \geq \frac{3\epsilon_H^2}{8L_3},$$

where $\mathbb{E}_\zeta$ denotes the expectation over the Rademacher random variable $\zeta$. Furthermore, if we choose $\delta \leq \epsilon_H/(3\epsilon_H + 8L_2)$, it holds that

$$\mathbb{E}[f(\widehat{\mathbf{y}}) - f(\mathbf{x})] \leq -\frac{\epsilon_H^2}{8L_3},$$

where $\mathbb{E}$ is over all randomness of the algorithm, and the total runtime is $\widetilde{O}\big((L_1^2/\epsilon_H^2)\mathbb{T}_h\big)$ if ApproxNC-Stochastic adopts online Oja's algorithm, and $\widetilde{O}\big((L_1^2/\epsilon_H^2)\mathbb{T}_g\big)$ if ApproxNC-Stochastic adopts Neon2$^{\text{online}}$.

## 5 Local Minima Finding Algorithms and Theories

In this section, we present our main algorithms to find approximate local minima for nonconvex finite-sum and stochastic optimization problems, based on the negative curvature descent algorithms proposed in previous section.

## 5.1 Finite-Sum Setting

We first consider the finite-sum setting in (4.5). The main idea is to apply a variance reduced stochastic gradient method until we encounter some first-order stationary points. Then we apply Algorithm 1 to find a negative curvature direction to escape any non-degenerate saddle point . The proposed algorithm is displayed in Algorithm 3. We call it FLASH for short.

To explain our algorithm in detail, we use the nonconvex Stochastic Variance Reduced Gradient (SVRG) method (Reddi et al., 2016; Allen-Zhu and Hazan, 2016; Lei et al., 2017) to find an $\epsilon$-first-order stationary point of (4.5). Without loss of generality, here we utilize a more general variance reduction method, the Stochastically Controlled Stochastic Gradient (SCSG) method (Lei et al., 2017), which reduces to the SVRG algorithm in special cases. When the algorithm reaches some $\epsilon$-first-order stationary point that satisfies $\|\nabla f(\mathbf{x})\|_2 \leq \epsilon$, we then apply the negative curvature descent, i.e., Algorithm 1, with a large step size. According to Lemma 4.3, if $\lambda_{\min}(\nabla^2 f(\mathbf{x})) < -\epsilon_H$, Algorithm 1 will escape from the saddle point. Otherwise, it will return the current iterate $\mathbf{x}$, which is already an approximate local minimum of (4.5). It is worth noting that the high-level structure of our algorithm bears a similarity with the algorithms proposed in Yu et al. (2017), namely GOSE. However, their algorithms cannot exploit the the third-order smoothness of the objective function, and therefore the step size of negative curvature descent in their algorithms is smaller. Yet they showed that with appropriate step size, their algorithms are able to escape each saddle point in one step.

---

**Algorithm 3** **F**ast **L**ocal minim**A** finding with third-order **S**moot**H**ness (FLASH-FiniteSum)

1: **Input:** $f(\cdot)$, $\mathbf{x}_0$, $L_1$, $L_2$, $L_3$, $\delta$, $\epsilon$, $\epsilon_H$, $K$
2: Set $\eta = 1/(6L_1 n^{2/3})$
3: **for** $k = 1, 2, ..., K$
4:     $\mathbf{g}_k \leftarrow \nabla f(\mathbf{x}_{k-1})$
5:     **if** $\|\mathbf{g}_k\|_2 > \epsilon$
6:         generate $T_k \sim \text{Geom}(n/(n+b))$
7:         $\mathbf{y}_0^{(k)} \leftarrow \mathbf{x}_{k-1}$
8:         **for** $t = 1, ..., T_k$
9:             randomly pick $i$ from $\{1, 2, ..., n\}$
10:            $\mathbf{y}_t^{(k)} \leftarrow \mathbf{y}_{t-1}^{(k)} - \eta(\nabla f_i(\mathbf{y}_{t-1}^{(k)}) - \nabla f_i(\mathbf{y}_0^{(k)}) + \mathbf{g}_k)$
11:         **end for**
12:         $\mathbf{x}_k \leftarrow \mathbf{y}_{T_k}^{(k)}$
13:     **else**
14:         $\mathbf{x}_k \leftarrow \text{NCD3-FiniteSum}(f(\cdot), \mathbf{x}_{k-1}, L_1, L_2, L_3, \delta, \epsilon_H)$
15:         **if** $\mathbf{x}_k = \perp$
16:             **return** $\mathbf{x}_{k-1}$
17: **end for**

---

Now we present runtime complexity analysis of Algorithm 3.

**Theorem 5.1.** Let $f(\mathbf{x}) = 1/n \sum_{i=1}^{n} f_i(\mathbf{x})$ and each component function $f_i(\mathbf{x})$ is $L_1$-smooth, $L_2$-Hessian Lipschitz continuous, and the third derivative of $f(\mathbf{x})$ is $L_3$-Lipschitz. If Algorithm 3 adopts Shifted-and-Inverted power method as its NCD3-FiniteSum subroutine to compute the negative curvature, then Algorithm 3 finds an $(\epsilon, \epsilon_H)$-second-order stationary point with probability at least $2/3$ in runtime

$$\widetilde{O}\bigg(\bigg(\frac{L_1\Delta_f n^{2/3}}{\epsilon^2} + \frac{L_3\Delta_f n}{\epsilon_H^2}\bigg)\mathbb{T}_h + \bigg(\frac{L_1^{1/2}L_3\Delta_f n^{3/4}}{\epsilon_H^{5/2}}\bigg)\mathbb{T}_g\bigg).$$

If Algorithm 3 uses Neon2$^{\text{svrg}}$ as its NCD3-FiniteSum subroutine, then it finds an $(\epsilon, \epsilon_H)$-second-order stationary point with probability at least $2/3$ in runtime

$$\widetilde{O}\bigg(\bigg(\frac{L_1\Delta_f n^{2/3}}{\epsilon^2} + \frac{L_3\Delta_f n}{\epsilon_H^2} + \frac{L_1^{1/2}L_3\Delta_f n^{3/4}}{\epsilon_H^{5/2}}\bigg)\mathbb{T}_g\bigg).$$

**Remark 5.2.** Note that the runtime complexity in Theorem 5.1 holds with constant probability $2/3$. In practice, one can always repeatedly run Algorithm 3 for at most $\log(1/\delta)$ times to achieve a higher confidence result that holds with probability at least $1 - \delta$.

**Remark 5.3.** Theorem 5.1 suggests that the runtime complexity of Algorithm 3 is $\widetilde{O}(n^{2/3}\epsilon^{-2} + n\epsilon_H^{-2} + n^{3/4}\epsilon_H^{-5/2})$ to find an $(\epsilon, \epsilon_H)$-second-order stationary point. In stark contrast, the runtime complexity of the state-of-the-art first-order finite-sum local minima finding algorithms (Allen-Zhu and Li, 2017b; Yu et al., 2017) is $\widetilde{O}(n^{2/3}\epsilon^{-2} + n\epsilon_H^{-3} + n^{3/4}\epsilon_H^{-7/2})$. Apparently, the runtime complexity of Algorithm 3 is better in the second and third terms. In particular, when $n \gtrsim \epsilon^{-3/2}$ and $\epsilon_H = \sqrt{\epsilon}$, the runtime complexity of Algorithm 3 is $\widetilde{O}(n^{2/3}\epsilon^{-2})$ and the runtime complexity of the state-of-the-art (Agarwal et al., 2016; Allen-Zhu and Li, 2017b; Yu et al., 2017) is $\widetilde{O}(n\epsilon^{-3/2})$. Therefore, Algorithm 3 is by a factor of $\widetilde{O}(n^{1/3}\epsilon^{1/2})$ faster than the state-of-the-art algorithms.

**Remark 5.4.** Note that the better performance guarantee of our algorithm stems from the exploitation of the third-order smoothness. Also notice that the runtime complexity of Algorithm 3 matches that of nonconvex SVRG algorithm (Reddi et al., 2016; Allen-Zhu and Hazan, 2016; Lei et al., 2017), which only finds first-order stationary points without the third-order smoothness assumption.

### 5.2 General Stochastic Setting

Now we consider the general nonconvex stochastic problem in (1.1). To find the local minimum, we also use SCSG (Lei et al., 2017), which is the state-of-the-art stochastic optimization algorithm, to find a first-order stationary point and then apply Algorithm 2 to find a negative curvature direction to escape the saddle point. The proposed method is presented in Algorithm 4. This algorithm is similar to Algorithm 3 in the finite-sum setting. However, we use a subsampled stochastic gradient $\nabla f_\mathcal{S}(\mathbf{x})$ in the outer loop (Line 4) of Algorithm 4, which is defined as $\nabla f_\mathcal{S}(\mathbf{x}) = 1/|\mathcal{S}| \sum_{i \in \mathcal{S}} \nabla F(\mathbf{x}; \xi_i)$. Another difference lies in the subroutine algorithm we employ to compute the negative curvature direction.

As shown in Algorithm 4, we use subsampleed gradient to determine whether the iterate $\mathbf{x}_{k-1}$ is a first-order stationary point or not. Suppose the stochastic gradient $\nabla F(\mathbf{x}; \xi)$ satisfies the gradient sub-Gaussian condition and the batch size $|\mathcal{S}_k|$ is large enough, then if $\|\nabla f_{\mathcal{S}_k}(\mathbf{x}_{k-1})\|_2 > \epsilon/2$,

**Algorithm 4** **F**ast **L**ocal minim**A** finding with third-order **S**moot**H**ness (FLASH-Stochastic)
---
1: **Input:** $f(\cdot)$, $\mathbf{x}_0$, $L_1$, $L_2$, $L_3$, $\delta$, $\epsilon$, $\epsilon_H$, $K$
2: Set $B \leftarrow \widetilde{O}(\sigma^2/\epsilon^2)$, $\eta = 1/(6L_1 B^{2/3})$
3: **for** $k = 1, 2, ..., K$
4:   uniformly sample a batch $\mathcal{S}_k \sim \mathcal{D}$ with $|\mathcal{S}_k| = B$
5:   $\mathbf{g}_k \leftarrow \nabla f_{\mathcal{S}_k}(\mathbf{x}_{k-1})$
6:   **if** $\|\mathbf{g}_k\|_2 > \epsilon/2$
7:     generate $T_k \sim \text{Geom}(B/(B+1))$
8:     $\mathbf{y}_0^{(k)} \leftarrow \mathbf{x}_{k-1}$
9:     **for** $t = 1, ..., T_k$
10:      randomly pick $\xi_i$ from distribution $\mathcal{D}$
11:      $\mathbf{y}_t^{(k)} \leftarrow \mathbf{y}_{t-1}^{(k)} - \eta(\nabla F(\mathbf{y}_{t-1}^{(k)}; \xi_i) - \nabla F(\mathbf{y}_0^{(k)}; \xi_i) + \mathbf{g}_k)$
12:    **end for**
13:    $\mathbf{x}_k \leftarrow \mathbf{y}_{T_k}^{(k)}$
14:  **else**
15:    $\mathbf{x}_k \leftarrow \text{NCD3-Stochastic}(f(\cdot), \mathbf{x}_{k-1}, L_1, L_2, L_3, \delta, \epsilon_H)$
16:    **if** $\mathbf{x}_k = \perp$
17:      **return** $\mathbf{x}_{k-1}$
18: **end for**

$\|\nabla f(\mathbf{x}_{k-1})\|_2 > \epsilon/4$ holds with high probability. Similarly, when $\|\nabla f_{\mathcal{S}_k}(\mathbf{x}_{k-1})\|_2 \leq \epsilon/2$, then with high probability $\|\nabla f(\mathbf{x}_{k-1})\|_2 \leq \epsilon$ holds.

Note that each iteration of the outer loop in Algorithm 4 consists of two cases: (1) if the norm of subsampled gradient $\nabla f_{\mathcal{S}_k}(\mathbf{x}_{k-1})$ is small, then we run one subroutine NCD3-Stochastic, i.e., Algorithm 2; and (2) if the norm of subsampled gradient $\nabla f_{\mathcal{S}_k}(\mathbf{x}_{k-1})$ is large, then we run one epoch of SCSG algorithm. This design can reduce the number of negative curvature calculations. There are two major differences between Algorithm 4 and existing algorithms in Xu and Yang (2017); Allen-Zhu and Li (2017b); Yu et al. (2017): (1) the step size of negative curvature descent step in Algorithm 4 is larger and improves the total runtime complexity; and (2) the minibatch size in each epoch of SCSG in Algorithm 4 can be set to 1 instead of being related to the accuracy parameters $\epsilon$ and $\epsilon_H$, while the minibatch size in each epoch of SCSG in the existing algorithms Xu and Yang (2017); Allen-Zhu and Li (2017b); Yu et al. (2017) has to depend on $\epsilon$ and $\epsilon_H$. Next we present the following theorem which spells out the runtime complexity of Algorithm 4.

**Theorem 5.5.** Let $f(\mathbf{x}) = \mathbb{E}_{\xi \sim \mathcal{D}}[F(\mathbf{x}; \xi)]$. Suppose the third derivative of $f(\mathbf{x})$ is $L_3$-Lipschitz, and each stochastic function $F(\mathbf{x}; \xi)$ is $L_1$-smooth and $L_2$-Hessian Lipschitz continuous. Suppose that the stochastic gradient $\nabla F(\mathbf{x}; \xi)$ satisfies the gradient sub-Gaussian condition with parameter $\sigma$. Set batch size $B = \widetilde{O}(\sigma^2/\epsilon^2)$. If Algorithm 4 adopts online Oja's algorithm as the NCD3-FiniteSum subroutine to compute the negative curvature, then Algorithm 4 finds an $(\epsilon, \epsilon_H)$-second-order stationary point with probability at least $1/3$ in runtime

$$\widetilde{O}\bigg(\bigg(\frac{L_1 \sigma^{4/3} \Delta_f}{\epsilon^{10/3}} + \frac{L_3 \sigma^2 \Delta_f}{\epsilon^2 \epsilon_H^2}\bigg)\mathbb{T}_h + \bigg(\frac{L_1^2 L_3 \Delta_f}{\epsilon_H^4}\bigg)\mathbb{T}_g\bigg).$$

If Algorithm 4 adopts Neon2$^{\text{online}}$ as the NCD3-FiniteSum subroutine, then it finds an $(\epsilon, \epsilon_H)$-second-order stationary point with probability at least $1/3$ in runtime

$$\widetilde{O}\bigg(\bigg(\frac{L_1\sigma^{4/3}\Delta_f}{\epsilon^{10/3}} + \frac{L_3\sigma^2\Delta_f}{\epsilon^2\epsilon_H^2} + \frac{L_1^2L_3\Delta_f}{\epsilon_H^4}\bigg)\mathbb{T}_g\bigg).$$

**Remark 5.6.** This runtime complexity in Theorem 5.5 again holds with a constant probability. In practice, one can repeatedly run Algorithm 4 for at most $\log(1/\delta)$ times to achieve a high probability result with probability at least $1-\delta$.

**Remark 5.7.** Theorem 5.5 suggests that the runtime complexity of Algorithm 4 is $\widetilde{O}(\epsilon^{-10/3} + \epsilon^{-2}\epsilon_H^{-2} + \epsilon_H^{-4})$ to find an $(\epsilon, \epsilon_H)$-second-order stationary point, compared with $\widetilde{O}(\epsilon^{-10/3} + \epsilon^{-2}\epsilon_H^{-3} + \epsilon_H^{-5})$ runtime complexity achieved by the state-of-the-art (Allen-Zhu and Li, 2017b; Yu et al., 2017). This indicates that the runtime complexity of Algorithm 4 is improved upon the state-of-the-art in the second and third terms. If we simply set $\epsilon_H = \sqrt{\epsilon}$, the runtime of Algorithm 4 is $\widetilde{O}(\epsilon^{-10/3})$ and that of the state-of-the-art stochastic local minima finding algorithms (Allen-Zhu, 2017; Tripuraneni et al., 2017; Xu and Yang, 2017; Allen-Zhu and Li, 2017b; Yu et al., 2017) becomes $\widetilde{O}(\epsilon^{-7/2})$, thus Algorithm 4 outperforms the state-of-the-art algorithms by a factor of $\widetilde{O}(\epsilon^{-1/6})$.

**Remark 5.8.** Note that we can set $\epsilon_H$ to a smaller value, i.e., $\epsilon_H = \epsilon^{2/3}$, and the total runtime complexity of Algorithm 4 remains $\widetilde{O}(\epsilon^{-10/3})$. It is also worth noting that the runtime complexity of Algorithm 4 matches that of the state-of-the-art stochastic optimization algorithm (SCSG) (Lei et al., 2017) which only finds first-order stationary points but does not impose the third-order smoothness assumption.

## 6 Conclusions

In this paper, we investigated the benefit of third-order smoothness of nonconvex objective functions in stochastic optimization. We illustrated that third-order smoothness can help faster escape saddle points, by proposing new negative curvature descent algorithms with improved decrement guarantees. Based on the proposed negative curvature descent algorithms, we further proposed practical stochastic optimization algorithms with improved run time complexity that find local minima for both finite-sum and stochastic nonconvex optimization problems.

## A Revisit of the SCSG Algorithm

In this section, for the purpose of self-containedness, we introduce the nonconvex stochastically controlled stochastic gradient (SCSG) algorithm (Lei et al., 2017) for general smooth nonconvex optimization problems with finite-sum structure, which is described in Algorithm 5.

The following lemma characterizes the function value gap after one epoch of Algorithm 5, which is a restatement of Theorem 3.1 in Lei et al. (2017).

**Lemma A.1.** (Lei et al., 2017) Let each $f_i(\cdot)$ be $L_1$-smooth. Set $\eta L_1 = \gamma(B/b)^{-2/3}$, $\gamma \le 1/6$ and

**Algorithm 5** SCSG $(f(\cdot), \mathbf{x}_0, T, \eta, B, b, \epsilon)$

1: **initialization:** $\widetilde{\mathbf{x}}_0 = \mathbf{x}_0$
2: **for** $k = 1, 2, ..., K$
3:     uniformly sample a batch $\mathcal{S}_k \subset [n]$ with $|\mathcal{S}_k| = B$
4:     $\mathbf{g}_k \leftarrow \nabla f_{\mathcal{S}_k}(\widetilde{\mathbf{x}}_{k-1})$
5:     $\mathbf{x}_0^{(k)} \leftarrow \widetilde{\mathbf{x}}_{k-1}$
6:     generate $T_k \sim \text{Geom}(B/(B+b))$
7:     **for** $t = 1, ..., T_k$
8:        randomly pick $\widetilde{\mathcal{I}}_{t-1} \subset [n]$ with $|\widetilde{\mathcal{I}}_{t-1}| = b$
9:        $\nu_{t-1}^{(k)} \leftarrow \nabla f_{\widetilde{\mathcal{I}}_{t-1}}(\mathbf{x}_{t-1}^{(k)}) - \nabla f_{\widetilde{\mathcal{I}}_{t-1}}(\mathbf{x}_0^{(k)}) + \mathbf{g}_k$
10:       $\mathbf{x}_t^{(k)} \leftarrow \mathbf{x}_{t-1}^{(k)} - \eta \nu_{t-1}^{(k)}$
11:     **end for**
12:     $\widetilde{\mathbf{x}}_k \leftarrow \mathbf{x}_{T_k}^{(k)}$
13: **end for**
14: **output:** Sample $\widetilde{\mathbf{x}}_K^*$ from $\{\widetilde{\mathbf{x}}_k\}_{k=1}^K$ uniformly.

$B \geq 9$. Then at the end of the $k$-th outer loop of Algorithm 5, it holds that

$$\mathbb{E}[\|\nabla f(\widetilde{\mathbf{x}}_k)\|_2^2] \leq \frac{5L_1}{\gamma}\left(\frac{b}{B}\right)^{1/3}\mathbb{E}[f(\mathbf{x}_0^{(k)}) - f(\widetilde{\mathbf{x}}_k)] + \frac{6\,\mathbb{1}\{B < n\}}{B}\mathcal{V},$$

where $\mathcal{V}$ is the upper bound on the variance of the stochastic gradient.

We will consider two cases of Lemma A.1. The first case is when $f(\cdot)$ has a finite-sum structure and we set $B = n$, which means we use the full gradient in line 4 of Algorithm 5. In this case the SCSG algorithm resembles nonconvex SVRG (Reddi et al., 2016; Allen-Zhu and Hazan, 2016). Specifically, we have the following corollary:

**Corollary A.2.** Let each $f_i(\cdot)$ be $L_1$-smooth. Set parameters $b = 1$, $B = n$, and $\eta = 1/(6L_1B^{2/3})$. Then at the end of the $k$-th outer loop of Algorithm 5, it holds that

$$\mathbb{E}[\|\nabla f(\widetilde{\mathbf{x}}_k)\|_2^2] \leq \frac{30L_1}{n^{1/3}}\mathbb{E}[f(\mathbf{x}_0^{(k)}) - f(\widetilde{\mathbf{x}}_k)].$$

The second case corresponds to the stochastic setting in (1.1). In this case, we have $\mathbf{g}_k = \nabla f_{\mathcal{S}_k}(\widetilde{\mathbf{x}}_{k-1}) = 1/B \sum_{i \in \mathcal{S}_k} \nabla F(\widetilde{\mathbf{x}}_{k-1}; \xi_i)$ and $n$ is relatively large, i.e., $n \gg O(1/\epsilon^2)$. Then we have the following corollary.

**Corollary A.3.** Let each stochastic function $F(\mathbf{x}; \xi)$ be $L_1$-smooth and suppose that $\nabla F(\mathbf{x}; \xi)$ satisfies the gradient sub-Gaussian condition in Definition 3.7. Suppose that $n \gg O(1/\epsilon^2)$ and $n > B$. Set parameters $b = 1$ and $\eta = 1/(6L_1B^{2/3})$. Then at the end of the $k$-th outer loop of Algorithm 5, it holds that

$$\mathbb{E}[\|\nabla f(\widetilde{\mathbf{x}}_k)\|_2^2] \leq \frac{30L_1}{B^{1/3}}\mathbb{E}[f(\widetilde{\mathbf{x}}_0^{(k)}) - f(\widetilde{\mathbf{x}}_k)] + \frac{12\sigma^2}{B}.$$

Note that by Vershynin (2010) the sub-Gaussian stochastic gradient in Definition 3.7 implies $\mathbb{E}[\|\nabla F(\mathbf{x};\xi) - \nabla f(\mathbf{x})\|] \leq 2\sigma^2$ and thus we take $\mathcal{V} = 2\sigma^2$ in Corollary A.3.

## B  Proofs for Negative Curvature Descent

In this section, we first prove the lemmas that characterizes the function value decrease in negative curvature descent algorithms.

### B.1  Proof of Lemma 4.3

*Proof.* Since by assumptions, $f(\mathbf{x})$ is $L_3$-Hessian Lipschitz continuous, according to Lemma 1 in Anandkumar and Ge (2016), for any $\mathbf{x}, \mathbf{y} \in \mathbb{R}^d$, we have

$$f(\mathbf{y}) \leq f(\mathbf{x}) + \langle \nabla f(\mathbf{x}), \mathbf{y} - \mathbf{x} \rangle + \frac{1}{2}(\mathbf{y} - \mathbf{x})^\top \nabla^2 f(\mathbf{x})(\mathbf{y} - \mathbf{x})$$
$$+ \frac{1}{6} \langle \nabla^3 f(\mathbf{x}), (\mathbf{y} - \mathbf{x})^{\otimes 3} \rangle + \frac{L_3}{24} \|\mathbf{y} - \mathbf{x}\|_2^4.$$

Denote the input point $\mathbf{x}$ of Algorithm 1 as $\mathbf{y}_0$. Suppose that $\widehat{\mathbf{v}} \neq \perp$. By Lemmas 4.1 and 4.2, the ApproxNC-FiniteSum($\cdot$) algorithm returns a unit vector $\widehat{\mathbf{v}}$ such that

$$\widehat{\mathbf{v}}^\top \nabla^2 f(\mathbf{y}_0) \widehat{\mathbf{v}} \leq -\frac{\epsilon_H}{2} \tag{B.1}$$

holds with probability at least $1 - \delta$ within $\widetilde{O}(n + n^{3/4}\sqrt{L_1/\epsilon_H})$ evaluations of Hessian-vector product or stochastic gradient. Define $\mathbf{u} = \mathbf{x} + \eta\widehat{\mathbf{v}}$ and $\mathbf{w} = \mathbf{x} - \eta\widehat{\mathbf{v}}$. Then it holds that

$$\langle \nabla f(\mathbf{y}_0), \mathbf{u} - \mathbf{y}_0 \rangle + \frac{1}{6}\langle \nabla^3 f(\mathbf{y}_0), (\mathbf{u} - \mathbf{y}_0)^{\otimes 3} \rangle + \langle \nabla f(\mathbf{y}_0), \mathbf{w} - \mathbf{y}_0 \rangle + \frac{1}{6}\langle \nabla^3 f(\mathbf{y}_0), (\mathbf{w} - \mathbf{y}_0)^{\otimes 3} \rangle = 0.$$

Furthermore, recall that we have $\widehat{\mathbf{y}} = \mathbf{x} + \zeta\eta\widehat{\mathbf{v}}$ in Algorithm 1 where $\zeta$ is a Rademacher random variable and thus we have $\mathbb{P}(\zeta = 1) = \mathbb{P}(\widehat{\mathbf{y}} = \mathbf{u}) = 1/2$ and $\mathbb{P}(\zeta = -1) = \mathbb{P}(\widehat{\mathbf{y}} = \mathbf{w}) = 1/2$, which immediately implies

$$\begin{aligned}
\mathbb{E}_\zeta[f(\widehat{\mathbf{y}}) - f(\mathbf{y}_0)] &\leq \frac{1}{2}(\widehat{\mathbf{y}} - \mathbf{y}_0)^\top \nabla^2 f(\mathbf{y}_0)(\widehat{\mathbf{y}} - \mathbf{y}_0) + \frac{L_3}{24}\|\widehat{\mathbf{y}} - \mathbf{y}_0\|_2^4 \\
&\leq \frac{\eta^2}{2}\widehat{\mathbf{v}}^\top \nabla^2 f(\mathbf{y}_0)\widehat{\mathbf{v}} + \frac{L_3\eta^4}{24}\|\widehat{\mathbf{v}}\|_2^4 \\
&\leq -\frac{\eta^2}{2}\frac{\epsilon_H}{2} + \frac{L_3\eta^4}{24}\|\widehat{\mathbf{v}}\|_2^4 \\
&= -\frac{3\epsilon_H^2}{8L_3}
\end{aligned} \tag{B.2}$$

holds with probability at least $1 - \delta$, where $\mathbb{E}_\zeta$ denotes the expectation over $\zeta$, the third inequality follows from (B.1) and in the last equality we used the fact that $\eta = \sqrt{3\epsilon_H/L_3}$. On the other hand, by $L_2$-smoothness of $f$ we have

$$\mathbb{E}_\zeta[f(\widehat{\mathbf{y}}) - f(\mathbf{y}_0)] \leq \frac{\eta^2}{2}\widehat{\mathbf{v}}^\top \nabla^2 f(\mathbf{y}_0)\widehat{\mathbf{v}} + \frac{L_3\eta^4}{24}\|\widehat{\mathbf{v}}\|_2^4$$

$$\leq \frac{\eta^2 L_2}{2} + \frac{L_3\eta^4}{24}\|\widehat{\mathbf{v}}\|_2^4$$

$$= \frac{3\epsilon_H(\epsilon_H + 4L_2)}{8L_3}. \tag{B.3}$$

Combining (B.2) and (B.3) yields

$$\mathbb{E}[f(\widehat{\mathbf{y}}) - f(\mathbf{y}_0)] \leq -\frac{3(1-\delta)\epsilon_H^2}{8L_3} + \frac{3\delta\epsilon_H(\epsilon_H + 4L_2)}{8L_3}$$

$$\leq -\frac{3(1-\delta)\epsilon_H^2}{16L_3},$$

where the second inequality holds if $\delta \leq \epsilon_H/(3\epsilon_H + 8L_2)$. Furthermore, plugging $\delta < 1/3$ into the above inequality we obtain $\mathbb{E}[f(\widehat{\mathbf{y}}) - f(\mathbf{y}_0)] \leq -\epsilon_H^2/(8L_3)$. □

### B.2 Proof of Lemma 4.6

*Proof.* The proof of Lemma 4.6 can be easily adapted from that of Lemma 4.3. The only difference lies in the runtime complexity since we employ ApproxNC-Stochastic($\cdot$) in Algorithm 2. By Lemmas 4.4 and 4.5, the total number of evaluations of Hessian-vector product or stochastic gradient of algorithm ApproxNC-Stochastic($\cdot$) is $\widetilde{O}(L_1^2/\epsilon_H^2)$. We omit the rest of the proof for simplicity. □

## C Proofs for Runtime Complexity of Algorithms

In this section, we prove the main theorems for our faster stochastic local minima finding algorithms.

### C.1 Proof of Theorem 5.1

*Proof.* We first compute the iteration complexity of the outer loop in Algorithm 3. Denote $\mathcal{I} = \{1, 2, \ldots, K\}$ as the index of iteration. We divide $\mathcal{I}$ into two disjoint categories $\mathcal{I}_1$ and $\mathcal{I}_2$. For $k \in \mathcal{I}_1$, we have that $\mathbf{x}_k$ is output by the NCD3-Stochastic phase of Algorithm 3, and for $k \in \mathcal{I}_2$, we have that $\mathbf{x}_k$ is output by the one-epoch SCSG phase of Algorithm 3. Clearly, $K = |\mathcal{I}| = |\mathcal{I}_1| + |\mathcal{I}_2|$, thus we will calculate $|\mathcal{I}_1|$ and $|\mathcal{I}_2|$ in sequence.

**Computing $|\mathcal{I}_1|$**: by Lemma 4.3 we have

$$\mathbb{E}[f(\mathbf{x}_{k-1}) - f(\mathbf{x}_k)] \geq \frac{\epsilon_H^2}{8L_3}, \quad \text{for } k \in \mathcal{I}_1. \tag{C.1}$$

Summing up (C.1) over $k \in \mathcal{I}_1$ yields

$$\frac{\epsilon_H^2|\mathcal{I}_1|}{8L_3} \leq \sum_{k \in \mathcal{I}_1}\mathbb{E}[f(\mathbf{x}_{k-1}) - f(\mathbf{x}_k)] \leq \sum_{k \in \mathcal{I}}\mathbb{E}[f(\mathbf{x}_{k-1}) - f(\mathbf{x}_k)] \leq \Delta_f,$$

where the second inequality use the following fact by Corollary A.2

$$0 \leq \mathbb{E}[\|\nabla f(\mathbf{x}_k)\|_2^2] \leq \frac{C_1 L_1}{n^{1/3}} \mathbb{E}[f(\mathbf{x}_{k-1}) - f(\mathbf{x}_k)], \quad \text{for } k \in \mathcal{I}_2, \tag{C.2}$$

and $C_1 = 30$ is a constant. Thus we further have

$$|\mathcal{I}_1| \leq \frac{8L_3 \Delta_f}{\epsilon_H^2}. \tag{C.3}$$

**Computing $|\mathcal{I}_2|$:** to get the upper bound of $|\mathcal{I}_2|$, we further decompose $\mathcal{I}_2$ as $\mathcal{I}_2 = \mathcal{I}_2^1 \cup \mathcal{I}_2^2$, where $\mathcal{I}_2^1 = \{k \in \mathcal{I}_2 \mid \|\mathbf{g}_k\|_2 > \epsilon\}$ and $\mathcal{I}_2^2 = \{k \in \mathcal{I}_2 \mid \|\mathbf{g}_k\|_2 \leq \epsilon\}$. Based on the definition, it is obvious that $\mathcal{I}_2^1 \cap \mathcal{I}_2^2 = \varnothing$ and $|\mathcal{I}_2| = |\mathcal{I}_2^1| + |\mathcal{I}_2^2|$. According to update rules of Algorithm 3, if $k \in \mathcal{I}_2^2$, then Algorithm 3 will execute the NCD3-FiniteSum stage in the $k$-th iteration, which indicates that $|\mathcal{I}_2^2| \leq |\mathcal{I}_1|$. Next we need to get the upper bound of $|\mathcal{I}_2^1|$. Summing up (C.2) over $k \in \mathcal{I}_2^1$ yields

$$\sum_{k \in \mathcal{I}_2^1} \mathbb{E}[\|\nabla f(\mathbf{x}_k)\|_2^2] \leq \frac{C_1 L_1}{n^{1/3}} \sum_{k \in \mathcal{I}_2^1} \mathbb{E}[f(\mathbf{x}_{k-1}) - f(\mathbf{x}_k)]$$

$$\leq \frac{C_1 L_1}{n^{1/3}} \sum_{k \in \mathcal{I}} \mathbb{E}[f(\mathbf{x}_{k-1}) - f(\mathbf{x}_k)]$$

$$\leq \frac{C_1 L_1 \Delta_f}{n^{1/3}},$$

where in the second inequality we use (C.1) and (C.2). Applying Markov's inequality we obtain with probability at least $2/3$

$$\sum_{k \in \mathcal{I}_2^1} \|\nabla f(\mathbf{x}_k)\|_2^2 \leq \frac{3 C_1 L_1 \Delta_f}{n^{1/3}}.$$

Recall that for $k \in \mathcal{I}_2^1$ in Algorithm 3, we have $\|\nabla f(\mathbf{x}_k)\|_2 > \epsilon$. We then obtain with probability at least $2/3$

$$|\mathcal{I}_2^1| \leq \frac{3 C_1 L_1 \Delta_f}{n^{1/3} \epsilon^2}. \tag{C.4}$$

**Computing Runtime:** To calculate the total computation complexity in runtime, we need to consider the runtime of stochastic gradient evaluations at each iteration of the algorithm. Specifically, by Lemma 4.3 we know that each call of the NCD3-FiniteSum algorithm takes $\widetilde{O}((n + n^{3/4}\sqrt{L_1/\epsilon_H})\mathbb{T}_h)$ runtime if FastPCA is used and $\widetilde{O}((n + n^{3/4}\sqrt{L_1/\epsilon_H})\mathbb{T}_g)$ runtime if Neon2 is used. On the other hand, Corollary A.2 shows that the length of one epoch of SCSG algorithm is $\widetilde{O}(n)$ which implies that the run time of one epoch of SCSG algorithm is $\widetilde{O}(n\mathbb{T}_g)$. Therefore, we can compute the total time complexity of Algorithm 3 with Shifted-and-Inverted power method as

follows

$$|\mathcal{I}_1| \cdot \widetilde{O}\bigg(\bigg(n + \frac{n^{3/4}L^{1/2}}{\epsilon_H^{1/2}}\bigg)\mathbb{T}_h\bigg) + |\mathcal{I}_2| \cdot \widetilde{O}(n\mathbb{T}_g)$$
$$= |\mathcal{I}_1| \cdot \widetilde{O}\bigg(\bigg(n + \frac{n^{3/4}L^{1/2}}{\epsilon_H^{1/2}}\bigg)\mathbb{T}_h\bigg) + (|\mathcal{I}_2^1| + |\mathcal{I}_2^2|) \cdot \widetilde{O}(n\mathbb{T}_g)$$
$$= |\mathcal{I}_1| \cdot \widetilde{O}\bigg(\bigg(n + \frac{n^{3/4}L^{1/2}}{\epsilon_H^{1/2}}\bigg)\mathbb{T}_h\bigg) + (|\mathcal{I}_2^1| + |\mathcal{I}_1|) \cdot \widetilde{O}(n\mathbb{T}_g),$$

where we used the fact that $|\mathcal{I}_2^2| \leq |\mathcal{I}_1|$ in the last equation. Combining the upper bounds of $|\mathcal{I}_1|$ in (C.3) and $|\mathcal{I}_2^1|$ in (C.4) yields the following runtime complexity of Algorithm 3 with Shifted-and-Inverted power method

$$O\bigg(\frac{L_3\Delta_f}{\epsilon_H^2}\bigg) \cdot \widetilde{O}\bigg(\bigg(n + \frac{n^{3/4}L^{1/2}}{\epsilon_H^{1/2}}\bigg)\mathbb{T}_h\bigg) + O\bigg(\frac{L_3\Delta_f}{\epsilon_H^2} + \frac{L_1\Delta_f}{n^{1/3}\epsilon^2}\bigg) \cdot \widetilde{O}(n\mathbb{T}_g)$$
$$= \widetilde{O}\bigg(\bigg(\frac{L_3\Delta_f n}{\epsilon_H^2} + \frac{L_1^{1/2}L_3\Delta_f n^{3/4}}{\epsilon_H^{5/2}}\bigg)\mathbb{T}_h + \bigg(\frac{L_1\Delta_f n^{2/3}}{\epsilon^2}\bigg)\mathbb{T}_g\bigg),$$

and the runtime of Algorithm 3 with Neon2$^{\text{svrg}}$ is

$$\widetilde{O}\bigg(\bigg(\frac{L_3\Delta_f n}{\epsilon_H^2} + \frac{L_1^{1/2}L_3\Delta_f n^{3/4}}{\epsilon_H^{5/2}} + \frac{L_1\Delta_f n^{2/3}}{\epsilon^2}\bigg)\mathbb{T}_g\bigg),$$

which concludes our proof. □

### C.2 Proof of Theorem 5.5

Before we prove theoretical results for the stochastic setting, we lay down the following useful lemma which states the concentration bound for sub-Gaussian random vectors.

**Lemma C.1.** (Ghadimi et al., 2016) Suppose the stochastic gradient $\nabla F(\mathbf{x}; \xi)$ is sub-Gaussian with parameter $\sigma$. Let $\nabla f_\mathcal{S}(\mathbf{x}) = 1/|\mathcal{S}| \sum_{i \in \mathcal{S}} \nabla F(\mathbf{x}; \xi_i)$ be a subsampled gradient of $f$. If the sample size $|\mathcal{S}| = 2\sigma^2/\epsilon^2(1 + \sqrt{\log(1/\delta)})^2$, then with probability $1 - \delta$,

$$\|\nabla f_\mathcal{S}(\mathbf{x}) - \nabla f(\mathbf{x})\|_2 \leq \epsilon$$

holds for any $\mathbf{x} \in \mathbb{R}^d$.

*Proof of Theorem 5.5.* Similar to the proof of Theorem 5.1, we first calculate the outer loop iteration complexity of Algorithm 4. Let $\mathcal{I} = \{1, 2, \ldots, K\}$ be the index of iteration. We use $\mathcal{I}_1$ and $\mathcal{I}_2$ to denote the index set of iterates that are output by the NCD3-FiniteSum stage and SCSG stage of Algorithm 4 respectively. It holds that $K = |\mathcal{I}| = |\mathcal{I}_1| + |\mathcal{I}_2|$. We will calculate $|\mathcal{I}_1|$ and $|\mathcal{I}_2|$ in sequence.

**Computing $|\mathcal{I}_1|$**: note that $|\mathcal{I}_1|$ is the number of iterations that Algorithm 4 calls NCD3-Stochastic to find the negative curvature. Recall the result in Lemma 4.6, for $k \in \mathcal{I}_1$, one execution of the

NCD3-Stochastic stage can reduce the function value up to

$$\mathbb{E}[f(\mathbf{x}_{k-1}) - f(\mathbf{x}_k)] \geq \frac{\epsilon_H^2}{8L_3}. \tag{C.5}$$

To get the upper bound of $|\mathcal{I}_1|$, we also need to consider iterates output by One-epoch SCSG. By Lemma A.3 it holds that

$$\mathbb{E}[\|\nabla f(\mathbf{x}_k)\|_2^2] \leq \frac{C_1 L_1}{B^{1/3}} \mathbb{E}[f(\mathbf{x}_{k-1}) - f(\mathbf{x}_k)] + \frac{C_2 \sigma^2}{B}, \quad \text{for } k \in \mathcal{I}_2, \tag{C.6}$$

where $C_1 = 30$ and $C_2 = 12$ are absolute constants. For $k \in \mathcal{I}_2$, we further decompose $\mathcal{I}_2$ as $\mathcal{I}_2 = \mathcal{I}_2^1 \cup \mathcal{I}_2^2$, where $\mathcal{I}_2^1 = \{k \in \mathcal{I}_2 \,|\, \|\mathbf{g}_k\|_2 > \epsilon/2\}$ and $\mathcal{I}_2^2 = \{k \in \mathcal{I}_2 \,|\, \|\mathbf{g}_k\|_2 \leq \epsilon/2\}$. It is easy to see that $\mathcal{I}_2^1 \cap \mathcal{I}_2^2 = \varnothing$ and $|\mathcal{I}_2| = |\mathcal{I}_2^1| + |\mathcal{I}_2^2|$. In addition, according to the concentration result on $\mathbf{g}_k$ and $\nabla f(\mathbf{x}_k)$ in Lemma C.1, if the sample size $B$ satisfies $B = O(\sigma^2/\epsilon^2 \log(1/\delta_0))$, then for any $k \in \mathcal{I}_2^1$, $\|\nabla f(\mathbf{x}_k)\|_2 > \epsilon/4$ holds with probability at least $1 - \delta_0$. For any $k \in \mathcal{I}_2^2$, $\|\nabla f(\mathbf{x}_k)\|_2 \leq \epsilon$ holds with probability at least $1 - \delta_0$. According to (C.6), we can derive that for any $k \in \mathcal{I}_2^1$,

$$\mathbb{E}[f(\mathbf{x}_{k-1}) - f(\mathbf{x}_k)] \geq \frac{B^{1/3}}{C_1 L_1} \mathbb{E}[\|\nabla f(\mathbf{x}_k)\|_2^2] - \frac{C_2 \sigma^2}{C_1 L_1 B^{2/3}}, \quad \text{for } k \in \mathcal{I}_2^1. \tag{C.7}$$

As for $|\mathcal{I}_2^2|$, because for any $k \in \mathcal{I}_2^2$, $\|\mathbf{g}_k\|_2 \leq \epsilon/2$, which will lead the algorithm to execute one step of NCD3-Stochastic stage in the next iteration, i.e., $k$-th iteration. Thus it immediately implies that $|\mathcal{I}_2^2| \leq |\mathcal{I}_1|$, and according to (C.6), we can also derive that

$$\mathbb{E}[f(\mathbf{x}_{k-1}) - f(\mathbf{x}_k)] \geq -\frac{C_2 \sigma^2}{C_1 L_1 B^{2/3}}, \quad \text{for } k \in \mathcal{I}_2^2. \tag{C.8}$$

Summing up (C.5) over $k \in \mathcal{I}_1$, (C.7) over $k \in \mathcal{I}_2^1$, (C.8) over $k \in \mathcal{I}_2^2$ and combining the results yields

$$\sum_{k \in \mathcal{I}} \mathbb{E}[f(\mathbf{x}_{k-1}) - f(\mathbf{x}_k)] \geq \sum_{k \in \mathcal{I}_1} \frac{\epsilon_H^2}{8L_3} + \frac{B^{1/3}}{C_1 L_1} \sum_{k \in \mathcal{I}_2^1} \mathbb{E}[\|\nabla f(\mathbf{x}_k)\|_2^2] - \sum_{k \in \mathcal{I}_2^1} \frac{C_2 \sigma^2}{C_1 L_1 B^{2/3}} - \sum_{k \in \mathcal{I}_2^2} \frac{C_2 \sigma^2}{C_1 L_1 B^{2/3}},$$

which immediately implies

$$\frac{|\mathcal{I}_1|\epsilon_H^2}{8L_3} + \frac{B^{1/3}}{C_1 L_1} \sum_{k \in \mathcal{I}_2^1} \mathbb{E}[\|\nabla f(\mathbf{x}_k)\|_2^2] \leq \Delta_f + \sum_{k \in \mathcal{I}_2^1} \frac{C_2 \sigma^2}{C_1 L_1 B^{2/3}} + \sum_{k \in \mathcal{I}_2^2} \frac{C_2 \sigma^2}{C_1 L_1 B^{2/3}}$$

$$\leq \Delta_f + |\mathcal{I}_2^1| \frac{C_2 \sigma^2}{C_1 L_1 B^{2/3}} + \frac{|\mathcal{I}_1| c_2 \sigma^2}{C_1 L_1 B^{2/3}},$$

where the first inequality uses the fact that $\Delta_f = f(\mathbf{x}_0) - \inf_{\mathbf{x}} f(\mathbf{x})$ and the second inequality is due to $|\mathcal{I}_2^2| \leq |\mathcal{I}_1|$. Applying Markov's inequality we have that

$$\frac{|\mathcal{I}_1|\epsilon_H^2}{8L_3} + \frac{B^{1/3}}{C_1L_1}\sum_{k\in\mathcal{I}_2^1}\|\nabla f(\mathbf{x}_k)\|_2^2 \leq 3\Delta_f + 3|\mathcal{I}_2^1|\frac{C_2\sigma^2}{C_1L_1B^{2/3}} + 3\frac{|\mathcal{I}_1|C_2\sigma^2}{C_1L_1B^{2/3}}$$

holds with probability at least $2/3$. Note that $\|\nabla f(\mathbf{x}_k)\|_2 \geq \epsilon/4$ with probability at least $1-\delta_0$. We conclude that by union bound we have that

$$\frac{|\mathcal{I}_1|\epsilon_H^2}{8L_3} + \frac{|\mathcal{I}_2^1|B^{1/3}\epsilon^2}{16C_1L_1} \leq 3\Delta_f + \frac{3|\mathcal{I}_2^1|C_2\sigma^2}{C_1L_1B^{2/3}} + \frac{3|\mathcal{I}_1|C_2\sigma^2}{C_1L_1B^{2/3}}$$

holds with probability at least $2/3(1-\delta_0)^{|\mathcal{I}_2^1|}$. We can set $B$ such that

$$\frac{3C_2\sigma^2}{C_1L_1B^{2/3}} \leq \frac{\epsilon_H^2}{16L_3},$$

which implies

$$B \geq \left(\frac{48C_2L_3\sigma^2}{C_1L_1}\right)^{3/2}\frac{1}{\epsilon_H^3}. \tag{C.9}$$

Combining the above two inequalities yields

$$\frac{|\mathcal{I}_1|\epsilon_H^2}{16L_3} + \frac{|\mathcal{I}_2^1|B^{1/3}\epsilon^2}{16C_1L_1} \leq 3\Delta_f + \frac{3|\mathcal{I}_2^1|C_2\sigma^2}{C_1L_1B^{2/3}} \tag{C.10}$$

holds with probability at least $2/3(1-\delta_0)^{|\mathcal{I}_2^1|}$. Therefore, it holds with probability at least $2/3(1-\delta_0)^{|\mathcal{I}_2^1|}$ that

$$|\mathcal{I}_1| \leq \frac{48L_3\Delta_f}{\epsilon_H^2} + \frac{48C_2L_3\sigma^2}{C_1L_1B^{2/3}\epsilon_H^2}|\mathcal{I}_2^1| = O\left(\frac{L_3\Delta_f}{\epsilon_H^2}\right) + \widetilde{O}\left(\frac{L_3\sigma^2}{L_1B^{2/3}\epsilon_H^2}\right)|\mathcal{I}_2^1|. \tag{C.11}$$

As we can see from the above inequality, the upper bound of $|\mathcal{I}_1|$ is related to the upper bound of $|\mathcal{I}_2^1|$. We will derive the upper bound on $|\mathcal{I}_2^1|$ later.

**Computing $|\mathcal{I}_2|$**: we have shown that $|\mathcal{I}_2^2| \leq |\mathcal{I}_1|$. Thus we only need to compute the cardinality of subset $\mathcal{I}_2^1 \subset \mathcal{I}_2$, where $\|\mathbf{g}_k\|_2 > \epsilon/2$ for any $k \in \mathcal{I}_2^1$. By Lemma C.1 we can derive that with probability at least $1-\delta_0$, it holds that $\|\nabla f(\mathbf{x}_k)\|_2 > \epsilon/4$. According to (C.10), we have

$$\frac{|\mathcal{I}_2^1|B^{1/3}\epsilon^2}{16C_1L_1} \leq \frac{|\mathcal{I}_1|\epsilon_H^2}{16L_3} + \frac{|\mathcal{I}_2^1|B^{1/3}\epsilon^2}{16C_1L_1} \leq 3\Delta_f + \frac{3|\mathcal{I}_2^1|C_2\sigma^2}{C_1L_1B^{2/3}} \tag{C.12}$$

holds with probability at least $2/3(1-\delta_0)^{|\mathcal{I}_2^1|}$. Further ensure that $B$ satisfies

$$\frac{3C_2\sigma^2}{C_1L_1B^{2/3}} \leq \frac{B^{1/3}\epsilon^2}{32C_1L_1},$$

which implies

$$B \geq \frac{96C_2\sigma^2}{\epsilon^2}. \tag{C.13}$$

Finally we get the upper bound of $|\mathcal{I}_2^1|$,

$$|\mathcal{I}_2^1| \leq \frac{96C_1 L_1 \Delta_f}{B^{1/3}\epsilon^2} = \widetilde{O}\left(\frac{L_1 \Delta_f}{\sigma^{2/3}\epsilon^{4/3}}\right),$$

where in the equation we use the fact in (C.9) and (C.13) that $B = \widetilde{O}(\sigma^2/\epsilon^2)$ and that $\epsilon_H = \sqrt{\epsilon}$.

We then plug the upper bound of $|\mathcal{I}_2^1|$ into (C.11) to obtain the upper bound of $\mathcal{I}_1$. Note that $B = \widetilde{O}(\sigma^2/\epsilon^2)$. Then we have

$$\begin{aligned}
|\mathcal{I}_1| &\leq \frac{48 L_3 \Delta_f}{\epsilon_H^2} + \frac{48 C_2 L_3 \sigma^2}{C_1 L_1 B^{2/3} \epsilon_H^2}|\mathcal{I}_2^1| \\
&= \widetilde{O}\left(\frac{L_3 \Delta_f}{\epsilon_H^2}\right) + \widetilde{O}\left(\frac{L_3 \sigma^{2/3} \epsilon^{4/3}}{L_1 \epsilon_H^2}\right) \cdot \widetilde{O}\left(\frac{L_1 \Delta_f}{\sigma^{2/3}\epsilon^{4/3}}\right) \\
&= \widetilde{O}\left(\frac{L_3 \Delta_f}{\epsilon_H^2}\right)
\end{aligned}$$

holds with probability at least $2/3(1-\delta_0)^{|\mathcal{I}_2^1|}$. Choosing sufficient small $\delta_0$ such that $(1-\delta_0)^{|\mathcal{I}_2^1|} > 1/2$, the upper bound of $\mathcal{I}_1$ and $\mathcal{I}_2^1$ holds with probability at least $1/3$.

**Computing Runtime**: By Lemma 4.6 we know that each call of the NCD3-Stochastic algorithm takes $\widetilde{O}((L_1^2/\epsilon_H^2)\mathbb{T}_h)$ runtime if FastPCA is used and $\widetilde{O}((L_1^2/\epsilon_H^2)\mathbb{T}_g)$ runtime if Neon2 is used. On the other hand, Corollary A.3 shows that the length of one epoch of SCSG algorithm is $\widetilde{O}(\sigma^2/\epsilon^2)$ which implies that the run time of one epoch of SCSG algorithm is $\widetilde{O}((\sigma^2/\epsilon^2)\mathbb{T}_g)$. Therefore, we can compute the total time complexity of Algorithm 3 with online Oja's algorithm as follows

$$\begin{aligned}
&|\mathcal{I}_1| \cdot \widetilde{O}\left(\frac{L_1^2}{\epsilon_H^2}\mathbb{T}_h\right) + |\mathcal{I}_2| \cdot \widetilde{O}\left(\frac{\sigma^2}{\epsilon^2}\mathbb{T}_g\right) \\
&= |\mathcal{I}_1| \cdot \widetilde{O}\left(\frac{L_1^2}{\epsilon_H^2}\mathbb{T}_h\right) + (|\mathcal{I}_2^1| + |\mathcal{I}_2^2|) \cdot \widetilde{O}\left(\frac{\sigma^2}{\epsilon^2}\mathbb{T}_g\right) \\
&= |\mathcal{I}_1| \cdot \widetilde{O}\left(\frac{L_1^2}{\epsilon_H^2}\mathbb{T}_h\right) + (|\mathcal{I}_2^1| + |\mathcal{I}_1|) \cdot \widetilde{O}\left(\frac{\sigma^2}{\epsilon^2}\mathbb{T}_g\right).
\end{aligned}$$

Plugging the upper bounds of $|\mathcal{I}_1|$ and $|\mathcal{I}_2^1|$ into the above equation yields the following runtime complexity of Algorithm 3 with online Oja's algorithm

$$\begin{aligned}
&\widetilde{O}\left(\frac{L_3 \Delta_f}{\epsilon_H^2}\right) \cdot \widetilde{O}\left(\frac{L_1^2}{\epsilon_H^2}\mathbb{T}_h\right) + \widetilde{O}\left(\frac{L_3 \Delta_f}{\epsilon_H^2} + \frac{L_1 \Delta_f}{\sigma^{2/3}\epsilon^{4/3}}\right) \cdot \widetilde{O}\left(\frac{\sigma^2}{\epsilon^2}\mathbb{T}_g\right) \\
&= \widetilde{O}\left(\left(\frac{L_1 \sigma^{4/3} \Delta_f}{\epsilon^{10/3}} + \frac{L_3 \sigma^2 \Delta_f}{\epsilon^2 \epsilon_H^2}\right)\mathbb{T}_h + \left(\frac{L_1^2 L_3 \Delta_f}{\epsilon_H^4}\right)\mathbb{T}_g\right),
\end{aligned}$$

and the runtime complexity of Algorithm 4 with Neon2$^{\text{online}}$ is

$$\widetilde{O}\bigg(\bigg(\frac{L_1\sigma^{4/3}\Delta_f}{\epsilon^{10/3}} + \frac{L_3\sigma^2\Delta_f}{\epsilon^2\epsilon_H^2} + \frac{L_1^2 L_3\Delta_f}{\epsilon_H^4}\bigg)\mathbb{T}_g\bigg),$$

which concludes our proof. □